\documentclass[12pt,a4paper,reqno]{amsart}

\usepackage{amsmath}
\usepackage{amssymb}
\usepackage{amsfonts}
\usepackage{xspace}

\advance\textheight by \topskip
\newtheorem{theorem}{Theorem}

\newtheorem{lemma}{Lemma}

\def\Prob{\mathbb{P}}

\def\bb#1{[\![#1]\!]}

\title[ Value and position of left--to--right
maxima]
{Combinatorics of geometrically distributed random
variables:\\  Value and position of the $r$th left--to--right
maximum}

\author{Arnold Knopfmacher and Helmut Prodinger}
\address{ Arnold Knopfmacher,
Centre for Applicable Analysis and Number Theory,
 Department of Applied Mathematics,
University of the Witwatersrand, P.~O. Wits, 
2050 Johannesburg, South Africa, email:
{\tt arnoldk@gauss.cam.wits.ac.za}.
}

\address{ Helmut Prodinger,
Centre for Applicable Analysis and Number Theory,
 Department of Mathematics,
University of the Witwatersrand, P.~O. Wits, 
2050 Johannesburg, South Africa, email:
{\tt helmut@gauss.cam.wits.ac.za}.
}

\date{December 21, 1998}

\begin{document}

\begin{abstract}
For words of length $n$, generated by independent geometric
random variables, we consider the average value and the
average position of the $r$th left--to--right maximum, for
fixed $r$ and $n\to \infty$.
\end{abstract}

\maketitle

\section{Introduction}
\label{sec:intro}

For a permutation $\sigma_1 \sigma_1\dots \sigma_n $,
a {\em left--to--right maximum\/} 
({\em outstanding element, record,\dots\/})
is an element $\sigma_j$ with $\sigma_j>\sigma_i$ 
for all $i=1,\dots,j-1$. The number of 
left--to--right maxima was first studied by R\'enyi \cite{Renyi62},
compare also \cite{Knuth68}. A survey of results on this topic
can be found in \cite{Glick78}.

Recently Wilf in \cite{Wilf95} proved the formula 
$(1-2^{-r})n$
for the average {\em value \/} of the $r$th 
left--to--right maximum, for fixed $r$ and $n\to\infty$;
for the average {\em position\/} he obtained
the asymptotic formula
$(\log n)^{r-1}/(r-1)!$.

In \cite{Prodinger96c} the number
of left--to--right maxima was investigated
in the model of {\em words\/} (strings) $a_1\dots a_n$, 
where the letters $a_i\in\mathbb N$ are independently
generated according to the geometric distribution with
$ \Prob{\{X=k\}}=pq^{k-1}$, with $p+q=1$.
(We find it useful also to use the abbreviation $Q=q^{-1}$.)
The motivation for this work came from Computer Science.
Also, since equal letters are now allowed, there are two
versions that should be considered in parallel, the standard
version, and the weak version, where `$<$' is replaced
by `$\le$,' which means that a new maximum  only
has to be larger or equal to the previous ones. 
The paper \cite{Prodinger96c} 
contains asymptotic results about the
average and the variance of the
 number
of left--to--right maxima in the context of geometric random
variables.
(H.--K.~Hwang and his collaborators obtained further
results about the limiting behaviour in \cite{BaHwLi98}.)

Motivated by Wilf's study we consider here the two 
parameters `value' and `position' of the 
$r$th left--to--right maximum
for geometric random variables. 
Summarizing our results, we obtain the 
asymptotic
formul{\ae}
$\frac rp$ and 
$\frac{1}{(r-1)!}\big(\frac pq\log_Qn\big)^{r-1}$
resp. 
$\frac {rq}p$ and 
$\frac{1}{(r-1)!}\big(p\log_Qn\big)^{r-1}$
in the weak case. 

A certain knowledge of \cite{Prodinger96c}
might be beneficial to understanding the present
derivations.

It should be noted that not all 
random strings of length $n$
{\em have\/} $r$
left--to--right maxima.

Let us start with the {\em value.} 
The generating function of interest is
\begin{equation*}
\frac{1}{1-z}
\prod_{i=1}^{h-1}\left\{
1+\frac{pq^{i-1}zu}{1-(1-q^i)z}\right\}
zpq^{h-1},
\end{equation*}
which originates from the (unique) decomposition of a string
as $a_1w_1\dots a_{r-1}w_{r-1}a_r w$ where
$a_1,\dots, a_r$ are the left--to--right maxima, 
the $w_i$ are the strings between them, and $w$ can be anything.
Note that if $a_k=l$, then this corresponds to a term
$pq^{l-1}zu$, and thus $w_k$ corresponds to $1/\big(
1-(1-q^k)z\big)$. A value $i$ must not necessarily occur
as a left--to--right maximum; that is reflected by the 
$1+\dots$ in the product. However, when we look for the
coefficient of $u^{r-1}$, we have seen $r-1$
 left--to--right maxima, and the $r$th has value $h$. 
What comes after that is irrelevant and covered by the
factor $1/(1-z)$.
(Compare \cite{Prodinger96c} for similar generating functions.)

In the sequel we find it useful
to use the abbreviation
$\bb{i}:=1-(1-q^i)z$.

The coefficients of $z^nu^{r-1}$, call them 
$\pi ^{(r)}_{n,h}$, are not probabilities, but 
${\pi ^{(r)}_{n,h}}\big/{\pi ^{(r)}_{n}}
$
are, 
 where $\pi ^{(r)}_{n}$ is the probability that a string of
length $n$ {\em has\/} 
$r$
left--to--right maxima; we find it as
\begin{equation*}
\pi ^{(r)}_{n}:=[z^nu^{r-1}]\frac{1}{1-z}
\sum_{h\ge1}\prod_{i=1}^{h-1}\left\{
1+\frac{pq^{i-1}zu}{\bb i}\right\}
zpq^{h-1}
=\sum_{h\ge1}\pi ^{(r)}_{n,h}.
\end{equation*}

Now we turn to the {\em position.}
Set
\begin{equation*}
\sigma ^{(r)}_{n,j}:=[z^nu^{r-1}v^j]\frac{1}{1-z}\sum_{h\ge1}
\prod_{k=1}^{h-1}\left\{
1+\frac{pq^{k-1}zvu}{1-(1-q^k)zv}\right\}
zvpq^{h-1},
\end{equation*}
then
${\sigma^{(r)}_{n,j}}\big/{\pi ^{(r)}_{n}}$
is the probability that a random string of length $n$ has
the $r$th maximum in position $j$. 
It is  the same decomposition as before, however, we are
not interested in the value $h$, so we sum over it. 
On the other hand, we label the position with the variable $v$, 
so we must make sure that every $z$ that does not appear in 
the factor $1/(1-z)$ must be multiplied by a $v$.
Computationally, we find it easier to work with the parameter
``position~$-r$,'' for which we have to consider
\begin{equation*}
\frac{1}{1-z}\sum_{h\ge1}
\prod_{k=1}^{h-1}\left\{
1+\frac{pq^{k-1}zu}{1-(1-q^k)zv}\right\}
zpq^{h-1},
\end{equation*}
since the variable $v$ appears in fewer places, as we don't have
to multiply all those $z$'s by $v$ which count for the $r$
 left--to--right maxima.

\section{Some technical lemmas}

In order to read off coefficients, we state the obvious
but nevertheless very useful formula

\begin{equation*}
[w^n]\sum_ia_if(b_iw)=\sum_ia_ib_i^n\cdot [w^n]f(w).
\end{equation*}
In all our applications, $\sum_ia_ib_i^n$ can be summed in
closed form.

\begin{lemma}
Assume that we have power series 
\begin{equation*}
A^{(j)}(w)=\sum_{n\ge1} a_{n}^{(j)}w^n, \qquad j=1,\dots,s.
\end{equation*}
Then
\begin{equation*}
[w^n]\sum_{1\le i_1<i_2<\dots<i_s}
A^{(1)}(wq^{i_1})\dots
A^{(s)}(wq^{i_s})=
\sum_{0=l_0<l_1<\dots<l_{s-1}<l_s=n}
\frac{a^{(s)}_{l_1-l_0}\dots a^{(1)}_{l_s-l_{s-1}}}
{(Q^{l_1}-1)\dots(Q^{l_s}-1)}.
\end{equation*}
\end{lemma}

\begin{proof}
For the sake of clarity, we treat the case $s=3$ and leave it
to the imagination of the reader to figure out the general case;
\begin{align*}
[w^n] &\sum_{1\le i<j<h}
A(wq^{i})B(wq^{j})
C(wq^{h})=
\sum_{l} \sum_{1\le i<j<h}
[w^{n-l}]A(wq^{i})B(wq^{j})\cdot[w^l]
C(wq^{h})\\
&=\sum_{l} \frac{1}{Q^l-1}
\sum_{1\le i<j}
[w^{n-l}]A(wq^{i})B(wq^{j})\cdot[w^l]
C(w)q^{jl}\\
&=\sum_{l} \frac{c_l}{Q^l-1}
\sum_{1\le i<j}
[w^{n}]A(wq^{i})B(wq^{j})
(wq^{j})^l\\
&=\sum_{l,m} \frac{c_l}{Q^l-1}
\sum_{1\le i<j}
[w^{n-m}]A(wq^{i})\cdot [w^m]B(wq^{j})
(wq^{j})^l\\
&=\sum_{l,m} \frac{c_l}{(Q^l-1)(Q^m-1)}
\sum_{1\le i}
[w^{n-m}]A(wq^{i})q^{im}\cdot [w^m]B(w)
w^l\\
&=\sum_{l,m} \frac{c_lb_{m-l}}{(Q^l-1)(Q^m-1)}
\sum_{1\le i}
[w^{n}]A(wq^{i})(wq^{i})^m\\
&=\sum_{l,m} \frac{c_l\,b_{m-l}\,a_{n-m}}{(Q^l-1)(Q^m-1)
(Q^n-1)}.
\end{align*}

\end{proof}

\begin{lemma}
\begin{align*}
[w^n]&\sum_{1\le i_1<i_2<\dots<i_s}
A^{(1)}(wq^{i_1})\dots
A^{(s)}(wq^{i_s})\,i_s\\
&=[t]
\sum_{0=l_0<l_1<\dots<l_{s-1}<l_s=n}
{a^{(s)}_{l_1-l_0}\dots a^{(1)}_{l_s-l_{s-1}}}
\prod_{i=1}^s
\bigg(\frac1{Q^{l_i}-1}
+t\frac{Q^{l_i}}{(Q^{l_i}-1)^2}
\bigg).
\end{align*}

\end{lemma}

\begin{proof}
The proof is essentially the same as before, if we note that
\begin{equation*}
\sum_{h>j}h\,q^{hl}=q^{jl}
\bigg(
 \frac{Q^l}{(Q^l-1)^2}+\frac{j}{Q^l-1}
\bigg).
\end{equation*}

\end{proof}

Our quantities will eventually come out as alternating sums, 
and the appropriate treatment of them is Rice's method
which is surveyed in \cite{FlSe95}; the key point
is the 
following Lemma.
\begin{lemma} 
Let $\mathcal C$ be a curve surrounding the points $1,2,\dots,n$ 
in the complex plane and let $f(z)$ be analytic inside $\mathcal C$. Then
$$
\sum_{k=1}^n \binom nk \, {(-1)}^k f(k)=
-\frac 1{2\pi i} \int_{\mathcal C} [n;z] f(z) dz, 
$$
where
$$
[n;z]=\frac{(-1)^{n-1} n!}{z(z-1)\dots(z-n)}=
\frac{\Gamma(n+1)\Gamma(-z)}{\Gamma(n+1-z)}.
$$

\end{lemma}

Extending the contour of integration it turns out that 
under suitable growth conditions on $f(z)$ (compare \cite{FlSe95})
the asymptotic expansion
of the alternating sum is given by
$$ 
\sum \text {Res} \big( [n;z] f(z)\big)+\text{smaller order terms}
$$
where the sum is taken over all poles $z_0$
different from $1,\dots,n$. 
Poles that lie more to the left lead to smaller terms in the
asymptotic expansion.

The range $1,\dots,n$ for the summation is not sacred; 
if we sum, for example, over $k=2,\dots,n$, the contour must
encircle $2,\dots,n$, etc.

\section{The probability that there are $r$ maxima}

Now we want to read off the $n$th coefficients of the power series
of interest. For this, it is beneficial to use the following formula:

\begin{equation*}
[z^n]f(z)=(-1)^n[w^n](1-w)^{n-1}f\Big(\frac{w}{w-1}\Big).
\end{equation*}
This form can be found in \cite{KiMaPr95} and is based on
ideas concerning the Euler transform in \cite{FlRi92}.
Then the quantities come out automatically as alternating
sums, and Rice's method can be applied.

\begin{align*}
\pi ^{(r)}_{n}&=[z^n]\frac{1}{1-z}
\sum_{1\le i_1< \dots <i_{r-1} <h}
\frac{pq^{i_{1}-1}z}{\bb{ i_{1}}}\dots
\frac{pq^{i_{r-1}-1}z}{\bb {i_{r-1}}}
zpq^{h-1}\\
&=
\Big(\frac{p}{q}\Big)^{r}[z^n]\frac{1}{1-z}
\sum_{1\le i_1< \dots <i_{r-1} <h}
\frac{q^{i_{1}}z}{\bb{ i_{1}}}\dots
\frac{q^{i_{r-1}}z}{\bb {i_{r-1}}}
zq^{h}\\
&=(-1)^r
\Big(\frac{p}{q}\Big)^{r}(-1)^n[w^n](1-w)^{n-1}
\sum_{1\le i_1< \dots <i_{r-1} <h}
\frac{q^{i_{1}}w\dots q^{i_{r-1}}w}{
(1-q^{i_{1}}w)\dots(1-q^{i_{r-1}}w)}
wq^{h}\\
&=(-1)^r
\Big(\frac{p}{q}\Big)^{r}
\sum_{k=r}^{n}\binom{n-1}{k-1}(-1)^k
[w^k]
\sum_{1\le i_1< \dots <i_{r-1} <h}
\frac{q^{i_{1}}w\dots q^{i_{r-1}}w}{
(1-q^{i_{1}}w)\dots(1-q^{i_{r-1}}w)}
wq^{h}.
\end{align*}
Now the evaluation of the inner sum can be done by our
Lemma; $k=n$, $r=s$, $A^{(1)}(w)=\dots=
A^{(r-1)}(w)=\frac{w}{1-w}$, $A^{(r)}(w)=w$.
Therefore
\begin{align*}
[w^k]&
\sum_{1\le i_1< \dots <i_{r-1} <h}
\frac{q^{i_{1}}w\dots q^{i_{r-1}}w}{
(1-q^{i_{1}}w)\dots(1-q^{i_{r-1}}w)}
wq^{h}\\
&=\sum_{0=l_0<1=l_1<l_2<\dots<l_r=k}
\frac{1}{(Q-1)(Q^{l_2}-1)\dots (Q^{l_r}-1)}\\
&=\frac qp\sum_{2\le l_2<\dots<l_r=k}
\frac{1}{(Q^{l_2}-1)\dots (Q^{l_r}-1)}.
\end{align*}
Thus
\begin{align*}
\pi ^{(r)}_{n}
=(-1)^{r-1}
\Big(\frac{p}{q}\Big)^{r-1}
\sum_{k=r-1}^{n-1}\binom{n-1}{k}(-1)^k
f(k)
\end{align*}
with

\begin{equation*}
f(k)=\sum_{2\le l_2<\dots<l_r=k+1}
\frac{1}{(Q^{l_2}-1)\dots (Q^{l_r}-1)}.
\end{equation*}

In order to apply Rice's method one needs the continuation 
of $f(k)$ to the complex plane.
Using {\em symmetric functions, } one can always 
represent such iterated summations by powersums
\cite{Knuth68}
\begin{equation*}
\vartheta (k):=\sum_{l= 2}^k\frac{1}{(Q^{l}-1)^d},
\end{equation*}
and the task is reduced to continue
this quantity $\vartheta (k)$
to the complex plane. For this, the standard way of doing
it is via
\begin{equation*}
\vartheta (z):=\sum_{l\ge 2}\frac{1}{(Q^{l}-1)^d}
-\sum_{l\ge1 }\frac{1}{(Q^{l+z}-1)^d}.
\end{equation*}

However, we only need the values
$f(0),\dots,f(r-2)$.

In \cite{GrKnPa94} we learn how such a sum has to be 
interpreted; we thus find $f(1)=\dots=f(r-2)=0$
and
\begin{equation*}
f(0)=(-1)^{r}\Big(\frac{1}{Q-1}\Big)^{r-1}.
\end{equation*}

Hence Rice's method and the pole at $z=0$
give us
\begin{equation*}
\pi ^{(r)}_{n}=1+O\Big(\frac1n\Big),
\end{equation*}
which is intuitively clear.

Note that there are poles at $z=-1+2\pi ik/\log Q$,
$k\in\mathbb Z$, 
and they lead to a periodic fluctuation of order $\frac 1n$;
this phenomenon is well--known and appears in many
places (compare \cite{FlSe95} and some other references).

\section{The average value of the $r$th maximum}

Now we can safely deal with the quantities
$\pi ^{(r)}_{n,h}$ alone, and  the so computated
average value $E^{(r)}_{n}$ will be correct within an error term
of the form $1+O(\frac1n)$.
We compute

\begin{align*}
 E^{(r)}_{n}&=[z^n]\frac{1}{1-z}
\sum_{1\le i_1< \dots <i_{r-1} <h}
\frac{pq^{i_{1}-1}z}{\bb{ i_{1}}}\dots
\frac{pq^{i_{r-1}-1}z}{\bb {i_{r-1}}}
zpq^{h-1}h\\
&=
\Big(\frac{p}{q}\Big)^{r}[z^n]\frac{1}{1-z}
\sum_{1\le i_1< \dots <i_{r-1} <h}
\frac{q^{i_{1}}z}{\bb{ i_{1}}}\dots
\frac{q^{i_{r-1}}z}{\bb {i_{r-1}}}
zq^{h}h\\
&=(-1)^r
\Big(\frac{p}{q}\Big)^{r}(-1)^n[w^n](1-w)^{n-1}
\sum_{1\le i_1< \dots <i_{r-1} <h}
\frac{q^{i_{1}}w\dots q^{i_{r-1}}w}{
(1-q^{i_{1}}w)\dots(1-q^{i_{r-1}}w)}
wq^{h}h\\
&=(-1)^r
\Big(\frac{p}{q}\Big)^{r}
\sum_{k=r}^{n}\binom{n-1}{k-1}(-1)^k
[w^k]
\sum_{1\le i_1< \dots <i_{r-1} <h}
\frac{q^{i_{1}}w\dots q^{i_{r-1}}w}{
(1-q^{i_{1}}w)\dots(1-q^{i_{r-1}}w)}
wq^{h}h.
\end{align*}
The evaluation of the inner sum is now done by the 
(second) Lemma:

\begin{align*}
[w^k]&
\sum_{1\le i_1< \dots <i_{r-1} <h}
\frac{q^{i_{1}}w\dots q^{i_{r-1}}w}{
(1-q^{i_{1}}w)\dots(1-q^{i_{r-1}}w)}
wq^{h}h\\
&=[t]\sum_{1=l_1<l_2<\dots<l_r=k}
\prod_{i=1}^r
\bigg(\frac1{Q^{l_i}-1}
+t\frac{Q^{l_i}}{(Q^{l_i}-1)^2}
\bigg).
\end{align*}
Or,

\begin{equation*}
E^{(r)}_n=
(-1)^{r-1}
\Big(\frac{p}{q}\Big)^{r}
\sum_{k=r-1}^{n-1}\binom{n-1}{k}(-1)^k
f(k)
\end{equation*}
with
\begin{equation*}
f(k)=[t]\sum_{1=l_1<l_2<\dots<l_r=k+1}
\prod_{i=1}^r
\bigg(\frac1{Q^{l_i}-1}
+t\frac{Q^{l_i}}{(Q^{l_i}-1)^2}
\bigg).
\end{equation*}
Again, $f(1)=\dots=f(r-2)=0$ and
\begin{align*}
f(0)&=[t]\sum_{1=l_1<l_2<\dots<l_r=1}
\prod_{i=1}^r
\bigg(\frac1{Q^{l_i}-1}
+t\frac{Q^{l_i}}{(Q^{l_i}-1)^2}
\bigg)\\&=(-1)^r
[t]\bigg(\frac1{Q-1}
+t\frac{Q}{(Q-1)^2}\bigg)^r=
(-1)^rr\frac{Q}{(Q-1)^{r+1}}
\end{align*}

Thus we have proved the following theorem

\begin{theorem}
The average value $E^{(r)}_n$ of the $r$th left--to--right maximum
in a random sequence of $n$ elements, generated by 
geometric random variables is given by
\begin{equation*}
E^{(r)}_n=\frac{r}{p}+O\Big(\frac1n\Big)
\qquad\text{for fixed $r$ and $n\to\infty.$}
\end{equation*}
\end{theorem}
A full asymptotic expansion would be available, at least in 
principle, with more involved computations, as well as the
variance.

Again, as in all the examples that will follow, the lower
order terms contain periodic fluctuations of the form
$\delta(\log_Qn)$.

\section{The average position of the $r$th maximum}

In order to compute this parameter (or rather the modified version),
we have to differentiate the generating function from the 
Introduction and plug in $v=1$.
The desired quantity is then obtained via
\begin{align*}
&[z^n]\Big(\frac pq\Big)^r
\frac{1}{1-z}[t]\sum_{1\le i_1<\dots<i_{r-1}<h}
\prod_{j=1}^{r-1}
\bigg(
\frac{q^{i_j}z}{\bb{i_j}}+t
\frac{q^{i_j}z(1-q^{i_j})z}{{\bb{i_j}}^2}
\bigg)zq^h\\
&=(-1)^n
[w^n](1-w)^{n-1}\Big(\frac pq\Big)^{r}(-1)^r
\times \\&\qquad\times 
[t]\sum_{1\le i_1<\dots<i_{r-1}<h}
\prod_{j=1}^{r-1}
\bigg(
\frac{q^{i_j}w}{1- q^{i_j}w}-t
\frac{q^{i_j}w(1-q^{i_j})w}{(1- q^{i_j}w)^2}
\bigg)wq^{h}\\
&=\Big(\frac pq\Big)^{r}(-1)^{r-1}
\sum_{k=r-1}^{n-1}\binom {n-1}k(-1)^kf(k)
\end{align*}

where
\begin{equation*}
f(k)=[w^{k+1}][t]\sum_{1\le i_1<\dots<i_{r-1}<h}
\prod_{j=1}^{r-1}
\bigg(
\frac{q^{i_j}w}{1- q^{i_j}w}-t
\frac{q^{i_j}w(1-q^{i_j})w}{(1- q^{i_j}w)^2}
\bigg)wq^{h}=f_1(k)+f_2(k)
\end{equation*}

with
\begin{equation*}
f_1(k)=[t][w^{k}]\sum_{1\le i_1<\dots<i_{r-1}<h}
\prod_{j=1}^{r-1}
\bigg(
\frac{q^{i_j}w}{1- q^{i_j}w}-t
\frac{q^{i_j}w}{(1- q^{i_j}w)^2}
\bigg)wq^{h}
\end{equation*}
and
\begin{equation*}
f_2(k)=[t][w^{k+1}]\sum_{1\le i_1<\dots<i_{r-1}<h}
\prod_{j=1}^{r-1}
\bigg(
\frac{q^{i_j}w}{1- q^{i_j}w}+t
\frac{(q^{i_j}w)^2}{(1- q^{i_j}w)^2}
\bigg)wq^{h}.
\end{equation*}
Now the two sums are in a form where our technical lemma 
applies!

For $f_1(k)$ note that $A_1(w)=\dots=A_{r-1}(w)=
\frac{w}{1-w}-\frac{tw}{(1-w)^2}$ and $A_r(w)=w$.
Thus
\begin{align*}
f_1(k)&=[t]
\sum_{1=l_1<l_2<\dots<l_r=k}\frac{
\big(1-t(l_2-l_1)\big)
\dots
\big(1-t(l_r-l_{r-1})\big)
}
{(Q^{l_1}-1)\dots(Q^{l_r}-1)}\\
&=-(k-1)
\sum_{1=l_1<l_2<\dots<l_r=k}\frac{1
}
{(Q^{l_1}-1)\dots(Q^{l_r}-1)}.
\end{align*}
For $f_2(k)$ note that $A_1(w)=\dots=A_{r-1}(w)=
\frac{w}{1-w}+\frac{tw^2}{(1-w)^2}$ and $A_r(w)=w$.
Thus											
\begin{align*}
f_2(k)&=[t]
\sum_{1=l_1<l_2<\dots<l_r=k+1}\frac{
\big(1+t(l_2-l_1-1)\big)
\dots
\big(1+t(l_r-l_{r-1}-1)\big)
}
{(Q^{l_1}-1)\dots(Q^{l_r}-1)}\\
&=(k-r)
\sum_{1=l_1<l_2<\dots<l_r=k+1}\frac{1
}
{(Q^{l_1}-1)\dots(Q^{l_r}-1)}.
\end{align*}
As we know from before, it is the ``value'' (the behaviour)
of $f(0)$ that is required. 
It is $f_1(0)$ that is dominant here:
Since in general
\begin{equation*}
\sum_{2\le l_2<\dots<l_{r-1}<0}a_{l_2}\dots
a_{l_{r-1}}=(-1)^r\frac{a_0^{r-1}-a_1^{r-1}}{a_0-a_1},
\end{equation*}

we find that as $z\to 0$
\begin{equation*}
f(z)\sim \frac{(-1)^r}{(Q-1)(Q^z-1)^{r-1}}.
\end{equation*}

Thus, according to the theory in \cite{FlSe95}, 
where it is explained in detail what kind of contribution an
$r$th order pole at $z=0$ gives,
we have proved that

\begin{theorem}

The average position of the $r$th left--to--right maximum
in a random sequence of $n$ elements, generated by 
geometric random variables is given by
\begin{equation*}
\frac{1}{(r-1)!}\Big(\frac pq\log_Qn\Big)^{r-1}
+O\big(\log^{r-2}n\big)
\qquad\text{for fixed $r$ and $n\to\infty.$}
\end{equation*}

\end{theorem}

\section{Weak left--to--right maxima; the value}

We mention here briefly the analogous developments for the 
instance of weak left--to--right maxima.

The generating function of interest is
\begin{equation*}
\frac{1}{1-z}
\prod_{i=1}^{h-1}
\bigg\{{1-\frac{pq^{i-1}zu}{\bb{i-1}}}\bigg\}^{-1}
pq^{h-1}z,
\end{equation*}
and the coefficient of $u^{r-1}$ therein is
\begin{align*}
&\frac{1}{1-z}\sum_{1\le i_1\le \dots\le i_{r-1}\le h}
\frac{pq^{i_{1}-1}z}{\bb{i_1-1}}\dots
\frac{pq^{i_{r-1}-1}z}{\bb{i_{r-1}}}pq^{h-1}z\\
&=p^r\frac{1}{1-z}\sum_{0\le i_1\le \dots\le i_{r-1}\le h}
\frac{q^{i_{1}}z}{\bb{i_1}}\dots
\frac{q^{i_{r-1}}z}{\bb{i_{r-1}}}q^{h}z.
\end{align*}
The technical lemmas that we need now are

\begin{lemma}
\begin{equation*}
[w^n]\sum_{0\le i_1\leq i_2\leq \dots \leq i_s}
A^{(1)}(wq^{i_1})\dots
A^{(s)}(wq^{i_s})=
\sum_{0=l_0<l_1<\dots<l_{s-1}<l_s=n}
\frac{a^{(s)}_{l_1-l_0}\dots a^{(1)}_{l_s-l_{s-1}}}
{(1-q^{l_1})\dots(1-q^{l_s})}
\end{equation*}
and
\begin{align*}
[w^n]&\sum_{0\le i_1 \leq i_2\leq \dots\leq i_s}
A^{(1)}(wq^{i_1})\dots
A^{(s)}(wq^{i_s})\,(i_s+1)\\
&=[t]
\sum_{0=l_0<l_1<\dots<l_{s-1}<l_s=n}
{a^{(s)}_{l_1-l_0}\dots a^{(1)}_{l_s-l_{s-1}}}
\prod_{i=1}^s
\bigg(\frac1{1-q^{l_i}}
+t\frac{q^{l_i}}{(1-q^{l_i})^2}
\bigg).
\end{align*}

\end{lemma}

We find

\begin{equation*}
E^{(r)}_n=
(-1)^{r-1}
p^{r}
\sum_{k=r-1}^{n-1}\binom{n-1}{k}(-1)^k
f(k)
\end{equation*}
with
\begin{equation*}
f(k)=[t]\sum_{1=l_1<l_2<\dots<l_r=k+1}
\prod_{i=1}^r
\bigg(\frac1{1-q^{l_i}}
+t\frac{q^{l_i}}{(1-q^{l_i})^2}
\bigg).
\end{equation*}

Also, $f(0)=(-1)^rr\frac{q}{(1-q)^{r+1}}$
and thus

\begin{theorem}

The average value $E^{(r)}_n$ of the $r$th left--to--right maximum
(in the weak sense)
in a random sequence of $n$ elements, generated by 
geometric random variables is given by
\begin{equation*}
E^{(r)}_n=\frac{r\,q}{p}+O\Big(\frac1n\Big)
\qquad\text{for fixed $r$ and $n\to\infty.$}
\end{equation*}

\end{theorem}

\section{Weak left--to--right maxima; the position}

The probability generating funtion of interest (up to normalization
by a factor that is basically $1$, as before)
is given by

\begin{equation*}
\frac{1}{1-z}\sum_{h\ge1}
\prod_{i=1}^{h-1}\left\{
1-\frac{pq^{i-1}zvu}{1-(1-q^i)zv}\right\}^{-1}
zvpq^{h-1},
\end{equation*}
and the desired expectated value is

\begin{align*}
[z^n] &p^r \frac{1}{1-z}
\sum_{0\leq i_1\leq \dots\leq i_{r-1}\leq h}
\prod_{j=1}^{r-1}
\bigg(
\frac{q^{i_j}z}{\bb{i_j}}+t
\frac{q^{i_j}z(1-q^{i_j})z}{{\bb{i_j}}^2}
\bigg)zq^h\\
&=p^r(-1)^{r-1}\sum_{k=r-1}^{n-1}\binom{n-1}{k}(-1)^kf(k)
\end{align*}
with $f(k)=f_1(k)+f_2(k)$ and

\begin{equation*}
f_1(k)=[t][w^{k}]\sum_{0\le i_1\le \dots\le i_{r-1}\le h}
\prod_{j=1}^{r-1}
\bigg(
\frac{q^{i_j}w}{1- q^{i_j}w}-t
\frac{q^{i_j}w}{(1- q^{i_j}w)^2}
\bigg)wq^{h}
\end{equation*}
and
\begin{equation*}
f_2(k)=[t][w^{k+1}]\sum_{0\le i_1\le \dots\le i_{r-1}\le h}
\prod_{j=1}^{r-1}
\bigg(
\frac{q^{i_j}w}{1- q^{i_j}w}+t
\frac{(q^{i_j}w)^2}{(1- q^{i_j}w)^2}
\bigg)wq^{h}.
\end{equation*}
We find
\begin{align*}
f_1(k)=-(k-1)
\sum_{1=l_1<l_2<\dots<l_r=k}\frac{1
}
{(1-q^{l_1})\dots(1-q^{l_r})}
\end{align*}
and
\begin{align*}
f_2(k)=(k-r)
\sum_{1=l_1<l_2<\dots<l_r=k+1}\frac{1
}
{(1-q^{l_1})\dots(1-q^{l_r})}.
\end{align*}

As $z\to 0$
\begin{equation*}
f(z)\sim \frac{(-1)^r}{(1-q)(1-q^z)^{r-1}},
\end{equation*}
and thus

\begin{theorem}

The average position of the $r$th left--to--right maximum
(in the weak sense)
in a random sequence of $n$ elements, generated by 
geometric random variables is given by
\begin{equation*}
\frac{1}{(r-1)!}\big(p\log_Qn\big)^{r-1}
+O\big(\log^{r-2}n\big)
\qquad\text{for fixed $r$ and $n\to\infty.$}
\end{equation*}

\end{theorem}

\bibliographystyle{plain}

\end{document}